\documentclass[12pt]{article}
\usepackage{amsmath,amssymb,comment}
\usepackage{amsfonts}
\usepackage{graphicx}
\newtheorem{theorem}{Theorem}[section]
\newtheorem{definition}[theorem]{Definition}
\newtheorem{remark}[theorem]{Remark}

\title{
Geometry of bifurcation sets: Exploring the parameter space
}

\begin{document}

\title{
Geometry of bifurcation sets: Exploring the parameter space
}

\author{R. Barrio\footnote{Departamento de Matem\'atica Aplicada and IUMA,
              Computational Dynamics group,
              University of Zaragoza,
              E-50009 Zaragoza, Spain, 
              rbarrio@unizar.es},
       S. Ib\'a\~{n}ez\footnote{Departamento de Matem\'aticas, University of
       Oviedo, E-33007 Oviedo, Spain, mesa@uniovi.es},
       and L. P\'erez\footnote{Departamento de Matem\'aticas, University of
       Oviedo, E-33007 Oviedo, Spain, perezplucia@uniovi.es}.}

\date{\today}
\maketitle

\begin{abstract}
Inspecting a $p$-dimensional parameter space by means of $(p-1)$-dimensional slices, changes can be detected that are only determined by the geometry of the manifolds that compose the bifurcation set. We refer to these changes as geometric bifurcations. They can be understood within the framework of the theory of singularities for differentiable mappings and, particularly, the Morse Theory. Working with a three-dimensional parameter space, geometric bifurcations are discussed in the context of two models of neuron activity: the Hindmarsh-Rose and the FitzHugh-Nagumo systems. Both are fast-slow systems with a small parameter that controls the time scale of a slow variable. Geometric bifurcations are observed on slices corresponding to fixed values of this distinguished small parameter.
\end{abstract}

\section{Introduction}
\label{sec:0}
Given a family of dynamical systems, the term bifurcation is used to refer to any qualitative change in the dynamics under variation of parameters. The stratification $B$ of the parameter space by bifurcation points is said the bifurcation set. On the other hand, when one explores a $q$-dimensional parameter space with $(q-s)$-parametric families of $s$-dimensional manifolds $S_\varepsilon$, where $s < q$ and $\varepsilon \in \mathbb{R}^{q-s}$, the geometry of $B$ itself can lead to changes in the sets $B_{\varepsilon}=B \cap S_{\varepsilon}$. We are interested in these changes, that we call \textit{geometric bifurcations}. As might be expected, the study of geometric bifurcations is closely related to the Theory of Singularities.

\begin{figure*}
\centerline{\includegraphics[width=1.\textwidth]{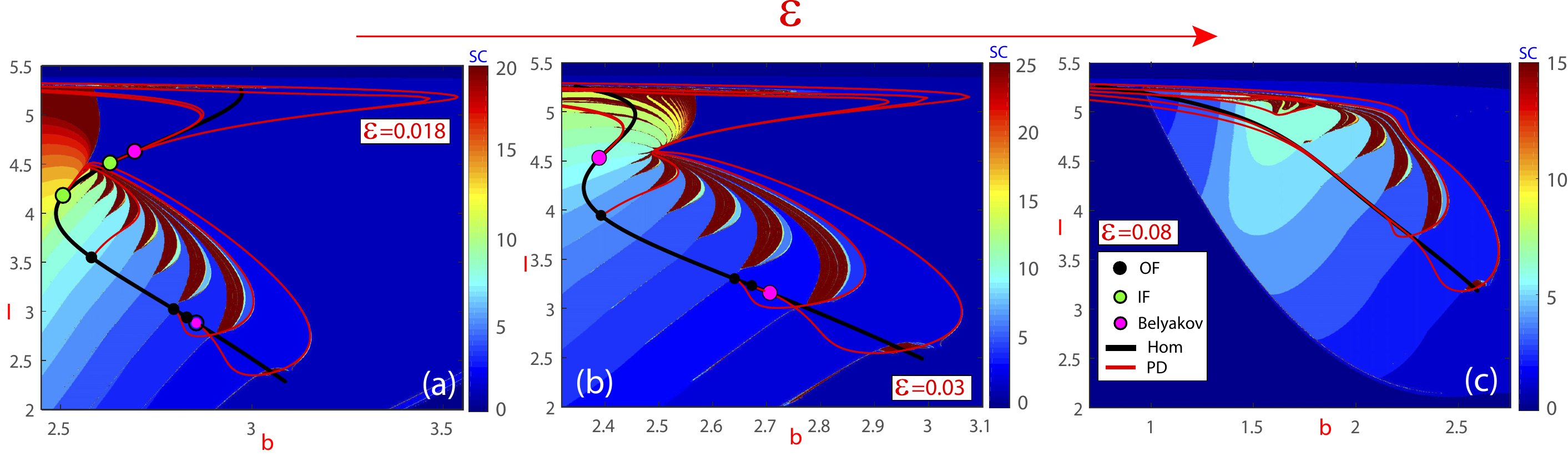}}
\caption{Spike-counting sweeping technique in the parameter plane $(b,I)$, with $\epsilon = 0.018$ (subplot (a)), $0.03$ (b) and $0.08$ (c). The color indicates the number of spikes per period, increasing from blue (0 spikes) to brown (the maximum value associated with chaotic behavior). Bifurcations: Homoclinic and period-doubling codimension-1 curves; orbit-flip (OF), inclination-flip (IF) and Belyakov codimension-2 points.}
\label{fig:general}
\end{figure*}

To illustrate the main aim of this article, the spike-counting (SC) technique (see \cite{barrio2011,Storace2008}) is used in Figure~\ref{fig:general} to show regions of periodic tonic spiking, chaotic bursting and regular bursting with different number of spikes for the classical Hindmarsh-Rose neuron model (more details in Subsection~\ref{ssA}). The SC sweeping method calculates for each set of parameter values and initial conditions the number of loops (spikes) of the stable limit cycle (when it exists), and this number is color-coded according to the number of spikes (vertical bar). Specifically, the dark blue regions correspond to stable spikings. The strips of a spike-adding cascade, ranging from blue to green, are related to fold/Hopf and fold/hom bursting (see \cite{BIPS2021} for details), which becomes chaotic in a chain of onion-like bulbs depicted in brown (see also \cite{BMSS14}). Additionally, several bifurcation lines are superimposed, such as the black one for homoclinic bifurcations and the red curves for period-doublings of periodic orbits. It is worth noting that the bifurcation curves of periodic orbits arise from homoclinic bifurcation points of codimension two, namely, the HR model exhibits orbit-flips (OF), inclination-flips (IF) and Belyakov bifurcations (see \cite{BIP2020} for a complete description). Furthermore, we highlight how the structure of the parameter plane $(b,I)$ changes as the small parameter $\varepsilon$ increases. An example would be the disappearance of some color bands and several bifurcation points of codimension two (e.g., the green and pink points) while, on the contrary, connections appear between bifurcation curves of periodic orbits (in red). Finally, we also look at how the Z-shape of the homoclinic curve (in black) evolves and disappears. Regarding these phenomena, the question that naturally arises is what is the mechanism underlying all these changes. As we will see later, the answer is that the bifurcation set goes through several geometric bifurcations as $\varepsilon$ varies.

Relevance of geometric bifurcations has been already pointed out in literature, although the terminology we are using in this paper is new. Geometric bifurcations explain the creation of isolas of homoclinic and heteroclinic bifurcations in \cite{algaba2003}. Namely, authors study Chua's equation and show how exploring the 3-parameter space with planar slices, a pair of codimension-two bifurcation points (T-points) collapse and disappear. When the collision occurs, the intersection between the slice and the codimension-two bifurcation curve is non-transversal. Through this process, spirals of codimension-one bifurcations approach together giving rise to isolas. In \cite{algaba2011}, a theoretical model explaining the behavior is provided. In \cite{algaba2015}, where authors study the distribution of homoclinic bifurcation curves in the Lorenz system, it is again evidenced how the evolution of a bifurcation set depends on the way in which the parameter space is explored.  In \cite{Wieczorek05,Wieczorek2007}, authors study a three-parameter family of three-dimensional vector fields that models an optically injected laser. Planar slices provide images of the bifurcation set where codimension-two bifurcations points collapse together and codimension-one bifurcation curves contact each other. Among other phenomena, isolas of homoclinic bifurcations are exhibited.

In this paper, we provide a minimal theoretical fra\-mework regarding geometric bifurcations by appealing to the theory of singularities on differentiable manifolds, and particularly, to Morse Theory \cite{katok08,matsumoto2002}. More importantly, it is shown how all the bifurcations that we introduce are observed in rather common neuron models. As already mentioned, we provide illustrations for the Hindmarsh-Rose \cite{HR84} and  FitzHugh-Nagumo systems \cite{Fitz61,Nagumo}. Both are fast-slow systems obtained as simplifications of the Hodgkin-Huxley model \cite{HH52} for the propagation of nerve impulses in axons. In both cases there are very rich bifurcation sets that, for instance, include homoclinic bifurcations of codimension one and two. In addition, we will see how geometric bifurcations arise when the parameter space is explored with slices corresponding to fixed values of a distinguished parameter $\varepsilon$, namely, the parameter that controls the time scale of the slow variable. These models illustrate a key fact, the study of geometric bifurcations only makes sense when working with families where there are distinguished parameters that determine a specific way of inspecting the parameter space; otherwise, either a reparametrization may remove singular points or tangencies can be avoided by choosing a different direction for sweepings.

In \cite{golsch1985}, Golubitsky and Schaeffer also use the term ``distinguished parameter'' to refer to a specific parameter whose variation can be seen as typical or natural in a given model. They work with scalar equations of the form
\begin{equation}
\label{escalar_equation}
G(x,\alpha,\beta)=0
\end{equation}
for a single unknown $x$ and with $\alpha\in\mathbb{R}$ and $\beta=(\beta_1,\ldots,\beta_{q-1})\in\mathbb{R}^{q-1}$. Typically, (\ref{escalar_equation}) stands for the steady-state equation of a certain dynamical problem. In this setting, they refer to $\alpha$ as a distinguished parameter and say that $\beta_1,\ldots,\beta_{q-1}$ are unfolding parameters. Golubitsky and Schaeffer wonder about how the sets
\begin{equation}
S_\beta=\{(x,\alpha)\in \mathbb{R}^2 \, : \, G(x,\alpha,\beta)=0\}
\end{equation}
change as $\beta$ varies. In this context, they characterize and classify the bifurcation points $(x_0,\alpha_0,\beta_0)$ ---singularities in their terminology--- where a change occurs. The notion of geometric bifurcation is also close to this theory of singularities. We consider families of dynamical systems depending on $q$ parameters $(\lambda,\varepsilon)\in \mathbb{R}^{q-s} \times \mathbb{R}^s$, with $s<q$ and explore the parameter space taking $(q-s)$-dimensional manifolds with $\varepsilon$ fixed. We wonder about how the intersection between these manifolds and the bifurcation set changes as $\varepsilon$ varies. In our setting, there is no state variable and $\varepsilon$ represents the distinguished parameter (or parameters in case that $s > 1$).
For the slow-fast models that we consider in this paper, the choice of distinguished parameter is indeed a ``natural'' one: the small parameter $\varepsilon$ that controls the behavior of the slow variable.
Note that the theory of singularities, as introduced in \cite{golsch1985}, allows to understand, for instance, the appearance of isola of equilibria in bifurcation diagrams (see \cite{isolas-1}).

The notion of geometric bifurcation incorporates a rather particular concept of codimension. Thus, bearing in mind a three-dimensional parameter space, we can find only three types of geometric bifurcations regarding codimension:

\begin{enumerate}
\item Codimension one-plus-one geometric bifurcations:

Those that occur on a surface of codimension-one bifurcations when the parameter space is explored with one-parametric families of two-dimensional ma\-nifolds.
\item Codimension two-plus-one geometric bifurcations:

Those that occur on a curve of codimension-two bifurcations when the parameter space is explored with one-parametric families of two-dimensional ma\-nifolds.
\item Codimension one-plus-two geometric bifurcations:

Those that occur on a surface of codimension-one bifurcations when the parameter space is explored with two-parametric families of one-dimensional ma-nifolds.
\end{enumerate}

Simple pictures that allow to illustrate the notion of geometric bifurcation are provided in Section \ref{sec:01}. They help to understand the particular way in which the codimension of a geometric bifurcation is specified (see details in \cite{Wieczorek05,Wieczorek2007}). The existence of a two-plus-one geometric bifurcation has subsidiary consequences on the codimension-one bifurcation sets that emerge from it, these effects are also discussed. In Section \ref{sec:1} we formalize some notions and provide a minimal theoretical support for dealing with geometric bifurcations. Specifically, one-plus-one and two-plus-one bifurcations are discussed. How geometric bifurcations arise in the aforementioned neural models is shown in Section \ref{sec:2}. Sketches in Section \ref{sec:01} and theoretical schemes provided in Section \ref{sec:1} should be useful to understand many of the phenomena exhibited in the neuron models. Finally, some conclusions are provided.

\section{Geometric bifurcations: illustrative examples}
\label{sec:01}

\begin{figure*}[tbh!]
\begin{center}
\includegraphics[width=0.8\textwidth]{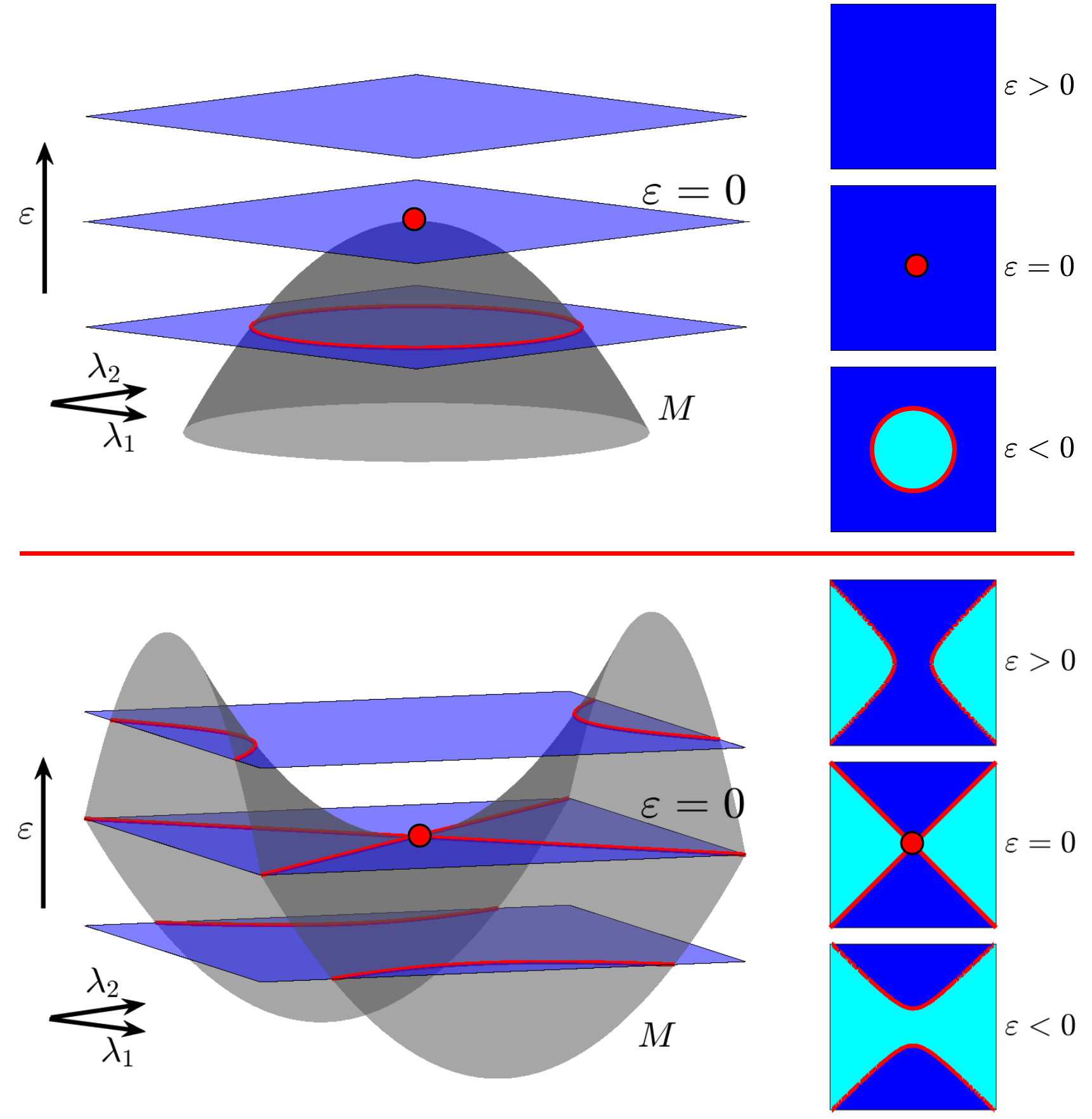}
\end{center}
\caption{One-plus-one geometric bifurcations exhibited in three-parameter spaces. Top: The bifurcation surface $M$ reaches a maximum value at $\varepsilon=0$. Exploring the surface with horizontal planes, an equivalence class (cyan) disappears because the limiting bifurcation curve (red circle) collapses to a point when $\varepsilon=0$ (top right panel). Bottom: The bifurcation surface has a saddle point. Passing through $\varepsilon=0$, the bifurcation curves (red branches of a hyperbola) change their position with respect to the asymptotes. As a consequence, the equivalence class around the central point switches (bottom right panel).}
\label{cod_1p1_3_parametros}
\end{figure*}

Let $X _{\lambda_1,\lambda_2,\varepsilon} $ be a three-parameter family of vector fields, and assume that a bifurcation surface $ M = \{(\lambda_1, \lambda_2, \varepsilon): \varepsilon+\lambda_1^2 + \lambda_2^2= 0 \} $ exists. When we explore the parameter space taking horizontal planes with $\varepsilon$ fixed,
we observe that, for $ \varepsilon <0$, there is a bifurcation curve, the circle $\lambda_1^2+\lambda_2^2=-\varepsilon$ (see top panel in Figure~\ref{cod_1p1_3_parametros}). In the absence of additional bifurcations, there is one equivalence class for values of the parameters inside the circle and one outside. When $\varepsilon=0$, the interior class disappears because the circle collapses to a point. Thus, when $\varepsilon>0$, bifurcations are no longer observed. If the bifurcation surface is given by
$M=\{(\lambda_1,\lambda_2,\varepsilon): \varepsilon-\lambda_1^2+\lambda_2^2 = 0 \}$, we observe a different behavior (see bottom panel in Figure~\ref{cod_1p1_3_parametros}). For
$\varepsilon<0$, there are two bifurcation curves, (the two branches of the hyperboloid $\varepsilon=\lambda_1^2-\lambda_2^2$). In the absence of  additional bifurcations, there exist two equivalence classes. When $\varepsilon=0$, the bifurcation set is given by $\lambda_1^2-\lambda_2^2 = 0$, two straight lines, and there is a singular point at $(0,0)$. For $\varepsilon> 0$, we find two curves again, but placed on different sectors. Note that in the passage through $\varepsilon=0$, the class which corresponds to the neighborhood of $(0,0)$ changes. In both cases, we say that there is a one-plus-one geometric bifurcation at $\varepsilon=0$. Note that intersections between the horizontal planes and the surface $M$ are all transverse except for $\varepsilon=0$.

\begin{figure*}[t]
\begin{center}
\includegraphics[width=0.48\textwidth]{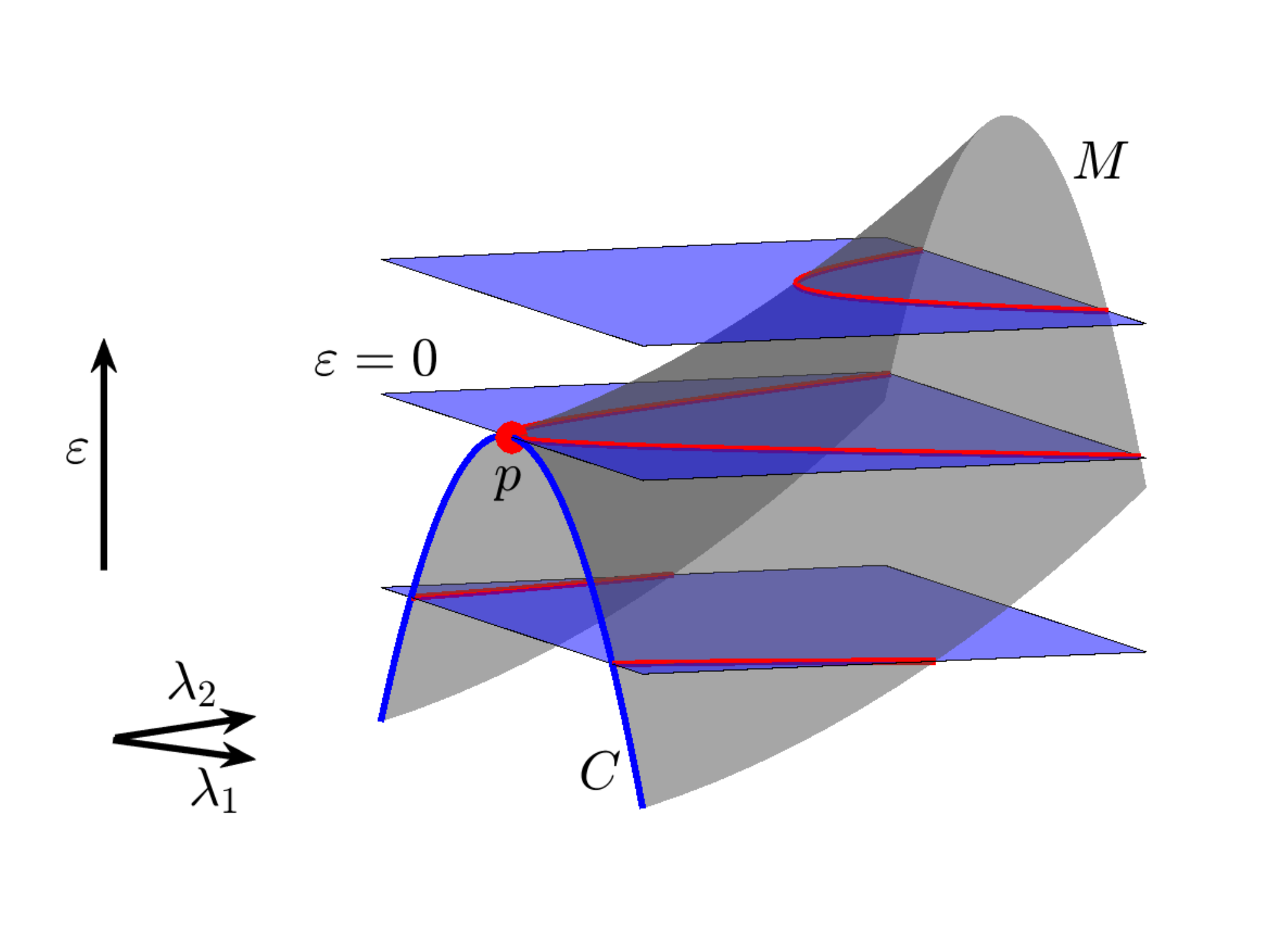}
\includegraphics[width=0.48\textwidth]{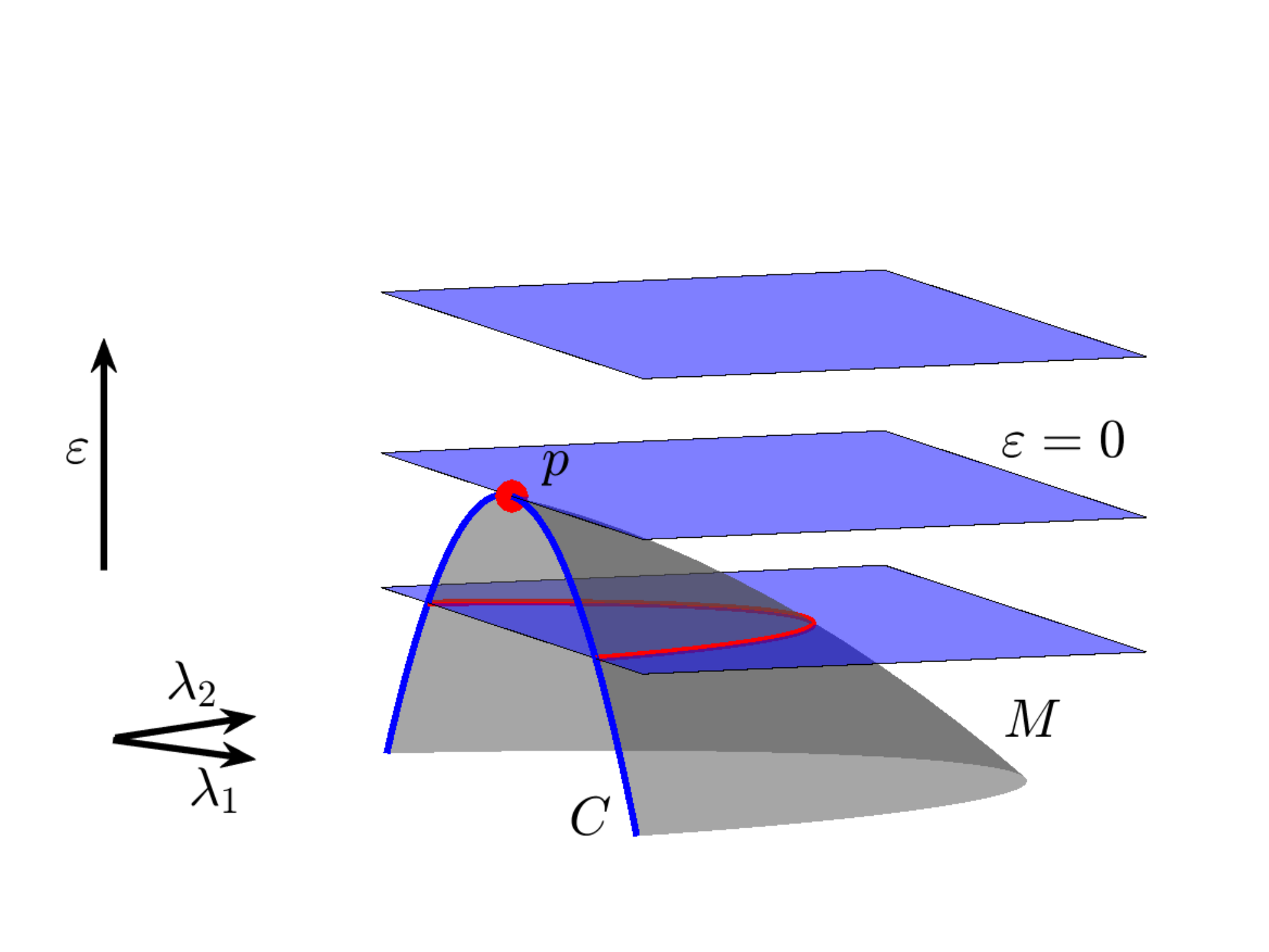}
\end{center}
\caption{Two-plus-one geometric bifurcations exhibited in three-parameter spaces and subsidiary effects. In both cases $C$ is a codimension-two bifurcation curve and $M$ is a codimension-one bifurcation surface wit a boundary along $C$. Left: Taking $\varepsilon$-slices two bifurcation curves (red) join in a unique one. Right: A bifurcation curve disappears after it collapses at the geometric bifurcation point.}
\label{cod_2p1_3_parametros}
\end{figure*}

In the previous examples, a geometric bifurcation appears because the three-parameter space is explored with planes and there is one which is tangent to a bifurcation surface. If there exists a codimension-two bifurcation curve exhibiting a point of tangency with one of planes, we say that there is a codimension--two-plus-one geometric bifurcation. Both panels in Figure \ref{cod_2p1_3_parametros} show a codimension-two bifurcation curve (blue) $C$. When the parameter space is explored with horizontal planes (fixed values of $\varepsilon$), there is a value of $\varepsilon$ for which the plane has a quadratic tangency with $C$ at a point $p$. Two codimension-two bifurcation points are exhibited for lower values of $\varepsilon$ and no one for bigger values. We say that there is a codimension--two-plus-one bifurcation at $p$.

More important, another interesting question is to wonder about the changes that occur in the bifurcation strata that arise from other strata of lower codimension. Both panels in Figure \ref{cod_2p1_3_parametros} show a surface $M$ (grey color) of codimension-one bifurcation points with a boundary given by the curve $C$. In the left panel, where the point $p$ is a saddle point on the surface, we can see how, as $\varepsilon$ increases, the horizontal slides show two curves (red color) that join together at $p$ to form a unique curve. On the other hand, looking at right panel we see that the point $p$ is a maximum on the bifurcation surface. In this case, as $\varepsilon$ increases, the horizontal slides show a unique curve (red color) that collapses at $p$ and disappear. The changes on the red curve are a subsidiary effect of the codimension two-plus-one geometric bifurcation. Note that in both cases, the bifurcation curves that are shown in each slide emerge from the codimension-two point with a well-defined tangent and we can observe how, ``generically'', they reconnect with each other or simply disappear, when the codimension-two points are no longer present. The former case is the one exhibited by the red curves in Figure \ref{fig:general} (see also examples in \cite{blackbeardthesis2012}). In Refs. \cite{algaba2003,algaba2011,Wieczorek05,Wieczorek2007} the case of spirals of bifurcations emerging at codimension-two bifurcation points was studied.

\begin{figure*}[t!]
\begin{center}
\includegraphics[width=1.\textwidth]{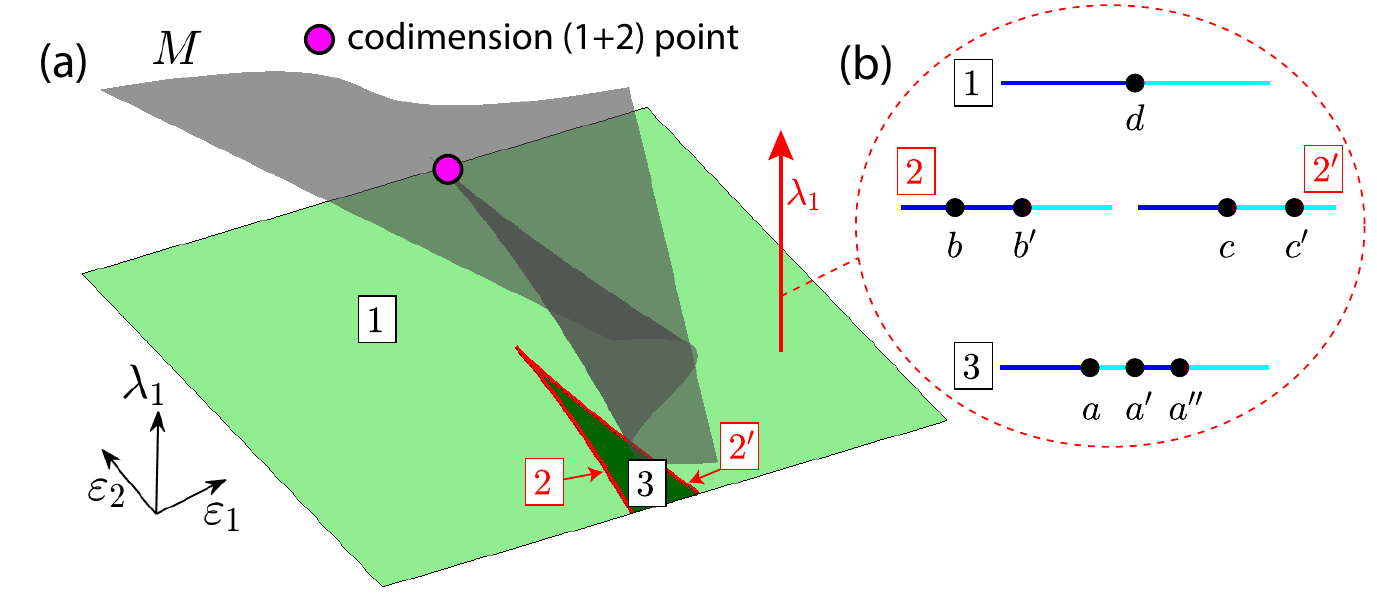}
\end{center}
\caption{(a) A one-plus-two geometric bifurcation. From the point $(\varepsilon_1,\varepsilon_2)=(0,0)$ two one-plus-one geometric bifurcation curves (red) emerge. They split the parameter space into two regions corresponding to non-equivalent bifurcation diagrams.
(b) Vertical one-dimensional ``slices'', obtained by fixing $\varepsilon_1$ and $\varepsilon_2$, are used to explore the bifurcation diagram exhibited in the panel (a), that is, we consider the three-parameter space as a two-parameter family of bifurcation diagrams in a one-parameter space. }
\label{cod_1p2}
\end{figure*}

New scenarios arise if we reduce the dimension of the slices.
Let $X_{\lambda_1,\varepsilon_1,\varepsilon_2}$ be a three-parameter family of vector fields and assume that the set $ M = \{(\lambda_1,\varepsilon_1,\varepsilon_2): \lambda_1^3+\lambda_1\varepsilon_2+\varepsilon_1=0\} $ is a bifurcation surface. Now, we explore the parameter space by taking vertical lines with $(\varepsilon_1,\varepsilon_2)$ fixed, see Figure~\ref{cod_1p2}. If $(\varepsilon_1,\varepsilon_2)$ belongs to the interior of region~3, i.e., the dark-green region bounded by the curve $27\varepsilon_1^2+4\varepsilon_2^3=0$ (in red), we find three bifurcation points $a$, $a'$ and $a''$ (Figure~\ref{cod_1p2}(b)---line~3). When parameters $(\varepsilon_1,\varepsilon_2)$ cross line~2 (Figure~\ref{cod_1p2}(b)), $a$ and $a'$ collapse at a codimension--one-plus-one bifurcation displayed in the two-parameter space we obtain by fixing $\varepsilon_2$. The same occurs along line $2'$, but the collision is between $a'$ and $a''$ instead. When parameters $(\varepsilon_1,\varepsilon_2)$ are in region~1 (Figure~\ref{cod_1p2}(b)---line~1), we find only one bifurcation point at $d$. We say that there is a codimension--one-plus-two geometric bifurcation at $(\varepsilon_1,\varepsilon_2)=(0,0)$. The Z-shape, exhibited by the bifurcation curves given by the intersection between $M$ and the vertical planes with $\varepsilon_2<0$ fixed, should be compared with the evolution of the black curves in Figure \ref{fig:general}. In addition, in \cite{CKKOS2007} it is illustrated the change in shape of homoclinic bifurcation curves, from Z-shape with two geometric folds to a curve without anyone (and so this is an example exhibiting a  one-plus-two geometric bifurcation in between), being relevant in fast-slow neuron systems.

\section{Geometric bifurcations and Morse Theory}
\label{sec:1}

Recall that a codimension-$k$ bifurcation is generic in the set of families of vector fields dependent on a number of parameters larger than or equal to $k$, but it can be avoided in families dependent on less parameters. A codimension-$k$ bifurcation is characterized by a set of $k$ degeneracy conditions independent of each other and a set of non-degeneracy (open) conditions. Any family satisfying these conditions is said a generic family or a generic unfolding.

Let $X_{\lambda,\varepsilon}:\mathbb{R}^n \rightarrow \mathbb{R}^n$ be a $\mathcal{C}^\infty$ family of vector fields with $(\lambda,\varepsilon)\in\mathbb{R}^q \times \mathbb{R}$. Assume that the family is a generic unfolding of a codimension-$k$ bifurcation at $(\lambda,\varepsilon)=(0,0)$, with $k\leq q$. Let $M$ be the $(q+1-k)$-dimensional smooth manifold through $(0,0)$ where such bifurcation is exhibited.
\begin{definition}
\label{geometric_bifurcation}
We say that the bifurcation point $(0,0)$ in the smooth manifold $M$ is \textbf{geometrically generic} with respect to $\varepsilon$ if the hyperplane $\varepsilon=0$ is transversal to $M$ at $(0,0)$. Otherwise, we say that the bifurcation point is \textbf{geometrically degenerate with respect to $\varepsilon$}.
\end{definition}

Let us recall the notion of Morse function (see \cite{Nicolaescu2011} for additional details).

\begin{definition}{\rm (\cite{Nicolaescu2011})}
Given a smooth manifold $M$ and a smooth function $f:M\rightarrow \mathbb{R}$, we say that a critical point $p_0$ of $f$ is \textbf{nondegenerate} if its Hessian is nondegenerate. We say that $f$ is a \textbf{Morse function} if all critical points are nondegenerate.
\end{definition}
Note that if $(0,0)$ is geometrically degenerate with respect to $\varepsilon$, then $(0,0)$ is a critical point of the height function $P_\varepsilon:M \rightarrow \mathbb{R}$ defined by $P_\varepsilon(\lambda,\varepsilon)=\varepsilon$ for each $(\lambda,\varepsilon)\in M$.
\begin{definition}
We say that the geometrically degenerate bifurcation at $(0,0)$ has \textbf{geometric codimension-one} with respect to $\varepsilon$ if the height function is locally a Morse function and $(0,0)$ is a critical point.
\end{definition}
The condition for geometric degeneracy is the fact that $(0,0)$ is a critical point of the height function $P_\varepsilon$. Imposing that $P_\varepsilon$ is a Morse function, we are setting a generic geometric condition, following \cite{Wieczorek05,Wieczorek2007}.

\begin{definition}
We say that a bifurcation point has \textbf{codimension $n$-plus-$1$} if it corresponds to a bifurcation of codimension-$n$ from a dynamical point of view, but the same point has codimension-$1$ from a geometrical point of view.
\end{definition}

\begin{remark}
To provide a theoretical framework for geometric bifurcations with codimension $m$, with $m>1$, one needs to consider $m$ distinguished parameters $(\varepsilon_1,\ldots,\varepsilon_m)$. In that context, the notion of \textbf{codimension $n$-plus-$m$} bifurcation would make sense. Nevertheless, in this paper, we do not enter into this formalization. However, note that in the previous section we have illustrated the case of a codimension-one-plus-two geometric bifurcation.
\end{remark}

In the sequel we assume that $q=2$. Since the height function is a Morse function locally around $(0,0)$, it follows from the Morse Lemma that choosing convenient coordinates $\lambda=(\lambda_1,\lambda_2)$, one can write either $
\varepsilon=\lambda_1^2+\lambda_2^2
$
or
$
\varepsilon=-\lambda_1^2-\lambda_2^2
$
or
$
\varepsilon=\lambda_1^2-\lambda_2^2.
$
In the first two cases, the level curves are circles (isolas) for $\varepsilon>0$ (resp. $\varepsilon<0$) and empty sets for $\varepsilon<0$ (resp. $\varepsilon>0$) (see Figure \ref{cod_1p1_3_parametros}(top)). We say isola-type to refer to this bifurcation. In the latter case the level curves correspond to a saddle case (see Figure \ref{cod_1p1_3_parametros}(bottom)) and we say simply saddle-type bifurcation. These two cases are also discussed in \cite{Wieczorek05,Wieczorek2007}.  The isola-type and saddle-type geometric bifurcations correspond to the codimension-one singularities named isola-center and simple bifurcation, respectively, in \cite{golsch1985}. Nevertheless, as we already mention in the introduction, there exist singularities which do not match any geometric bifurcation, the hysteresis point, a codimension-one singularity in \cite{golsch1985}, is one example.

Regarding codimension--two-plus-one bifurcations in a three-parameter space, in principle we only distinguish one type. Indeed, given a curve $C$ corresponding to a codimension-two bifurcation with a folding point, we can apply again the Morse Lemma and choose a convenient coordinate $\mu$ along $M$ such that $\varepsilon=\pm \mu^2$. Therefore, the level sets are given by two points that collapse and disappear (compare with the illustration provided in Figure \ref{cod_2p1_3_parametros}).

\begin{figure*}
\centerline{\includegraphics[width=0.9\textwidth]{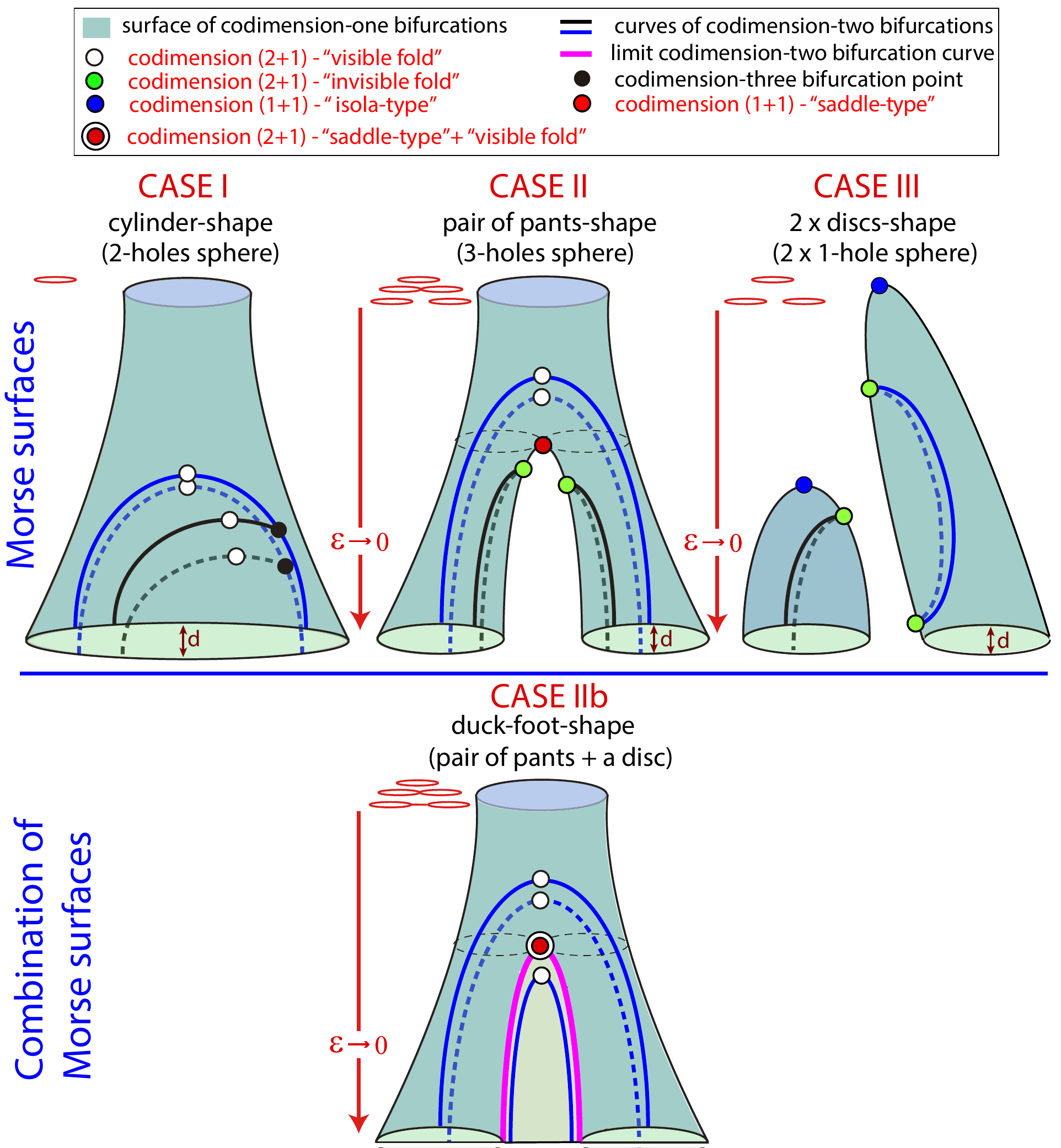}}
\caption{Theoretical scheme of the different topologies (Morse classification) of the bifurcation surfaces and the possible geometric and dynamic bifurcations on them. }
\label{topologia-def}
\end{figure*}

Another relevant point is to study how the different geometric bifurcations can appear in the bifurcation diagram of a given family of dynamical systems. We need to understand how a bifurcation surface exhibiting geometric bifurcations can look like. Figure \ref{topologia-def} shows a scheme of different scenarios. A key point is the topology of the bifurcation surface, having one or several tubular structures or even maxima. We classify the different 2D manifolds using the Morse theory \cite{matsumoto2002,Nicolaescu2011}. In all cases the value of a parameter $\varepsilon$ on the surface is our Morse function.

The different possibilities for Morse manifolds shown in Figure \ref{topologia-def} are illustrated later in Section~\ref{sec:2}, namely, Cases~I to III in Subsection~\ref{ssA} and Case~IIb in Subsection~\ref{ssB}. ---\emph{Case~I}---  surface is topologically equivalent to a cylinder, that is, a two-holes sphere. All values of $\varepsilon$ are regular and hence the structure of the surface is the same, a single circle, for all level sets. ---\emph{Case II}---  surface, a three-holes sphere, is called a ``pair of pants'' surface in the context of topology. In this case the Morse function given by the value of $\varepsilon$ has a saddle point (the red point in the middle plot of Figure~\ref{topologia-def}). This point is a codimension--one-plus-one geometric bifurcation point. In this case, the passage through the saddle corresponds to a connecting transition: two isolas collide in the saddle point and give rise to a unique isola. In the case shown in ---\emph{Case III}---, the pair of pants structure is broken, giving rise to two independent surfaces that disappear finally at a maximum of the surface when $\varepsilon$ grows (see the blue points). This case is a two discs shape (two one-hole spheres). Again the points of maxima are codimension--one-plus-one bifurcation points. Finally, ---\emph{Case IIb}---, that  we call ``duck-foot surface'', is obtained from the combination of a ``pair of pants'' and a ``disc'' (two Morse surfaces that meet along a curve). Note that this case requires the existence of a codimension-two bifurcation curve (magenta color) that acts as a limit of the tubular structures and enables the continuation in a disc surface. Therefore a cut on a transversal plane gives us the picture of two isolas connected with a curve. More combinations of Morse surfaces are possible, but we only discuss the one detected in the system of Subsection~\ref{ssB}.

Bifurcation surfaces include codimension-two bifurcation curves (black and blue lines) and codimension-three bifurcation points (black points). The value of $\varepsilon$ along the curves of codimension-two bifurcation plays a crucial role in the understanding of the global picture. Therefore, this parameter has a relevant ``physical'' meaning as it is assumed to be the parameter that can be changed in the system giving different phase space dynamics. All points along a curve of codimension-two bifurcation where $\varepsilon$ reaches an extremum (white and green points) are codimension--two-plus-one geometric bifurcations.

The tubular structures exhibited by the Morse surfaces allow us to distinguish between two different types of bifurcation curves of codimension two:  those that appear in pairs, on opposite sides of the cylinder; and those that appear as single curves with a fold point where they cross from one side of the cylinder to the other. All curves shown in Case I and those in blue color that we observed in the pair of pants of cases II and IIb are of the first type. The second type corresponds to the bifurcation curves in black color shown in Case II and all the curves in Case III. An independent case is the bifurcation curve (blue color) which is contained in the disk component of the Morse surface in  Case IIb. In all cases, the points on the curves where $\varepsilon$ reaches a maximum on the curves are codimension--two-plus-one geometric bifurcations. White points are folds with both branches contained in a unique side of the cylinder (or inside the disk for the Case IIb). Green points are folds where the two branches belong to different sides of the cylinder.

Let $d$ be the distance between both sides of a tubular structure, as depicted in Figure~\ref{topologia-def}. If $d$ is small enough, the two branches arising from green folds are indistinguishable and the folding point may be wrongly perceived as an end point. This situation happens in the fast-slow examples provided in the next section, where $d$ is exponentially small.  To distinguish these ``false'' end points, green points are said ``invisible folds'', in contrast with white points that are said ``visible folds'' (note that this situation is for $d\ll 1$). In Section \ref{sec:2}, we do not include models with $d$ big, but they would not provide any other noteworthy behavior either. Note that the appearance of invisible folds is linked to the formation of the ``pants'' or the splitting of a Morse surface in different components. This allows two disconnected isolas of a codimension-two bifurcation curve to meet and create a unique isola.

\section{Some models exhibiting geometric bifurcations}
\label{sec:2}

In this section we will consider two models that exhibit geometric bifurcations, namely, the Hindmarsh-Rose model and the FitzHugh-Nagumo system, showing how geometric bifurcations permit to visualize the global parameter space.

\subsection{The Hindmarsh-Rose model}
\label{ssA}
The Hindmarsh-Rose (HR) system was introduced in
\cite{HR84} as a reduction of the Hodgkin-Huxley equations \cite{HH52} to model the neuron behavior.
The HR model is described by three coupled nonlinear ODEs:

\begin{equation}
\label{HRmodel}
\left\{\begin{array}{l}
\dot{x}=y-ax^3+bx^2-z+I, \\
\dot{y}=c-dx^2-y, \\
\dot{z}=\varepsilon[s(x-x_0)-z],
\end{array}\right.
\end{equation}

\noindent where $x$ is the membrane potential, $y$ the fast and $z$ the slow gating variables for ionic current. Typically, $b$, $I$ and $\varepsilon$ are considered as free parameters, whereas the remaining ones are set as follow: $a=1$, $c=1$, $d=5$, $s=4$, $x_0=-1.6$.

The HR model has been deeply studied in recent years (see \cite{Barrio2017,BIP2020,BIPS2020,BIPS2021,BMSS14,barrio2011,Linaro2012,Shilnikov2008,Storace2008}) using different techniques. The system is particularly useful to understand the bursting phenomena and, in particular, the mechanisms of spike-adding. Behind these dynamics there are homoclinic bifurcations \cite{BMSS14,barrio2011,Linaro2012,Shilnikov2008,Storace2008}. In \cite{Barrio2017,BIP2020} we unveiled the global homoclinic structure and also part of the tangle of bifurcations of periodic orbits that emerge from the homoclinic skeleton, which allowed us to explain the spike-adding processes (see \cite{BIP2020,BIPS2021}).

In most of fast-slow systems with explicit small parameters, these parameters play a significant role and under their variation drastic changes in the global phase space are exhibited. The HR model is a paradigmatic example of this fact. As the small parameter $\varepsilon$ increases, numerous changes occur. In what follows, we show how the HR model exhibits geometric bifurcations with respect to $\varepsilon$ that, linked with dynamic bifurcations, explain these changes. Therefore, the HR model shows how geometric bifurcations permit us to help in the visualization of the global bifurcation diagram in the parameter space.

\begin{figure*}
\centerline{\includegraphics[width=1.\textwidth]{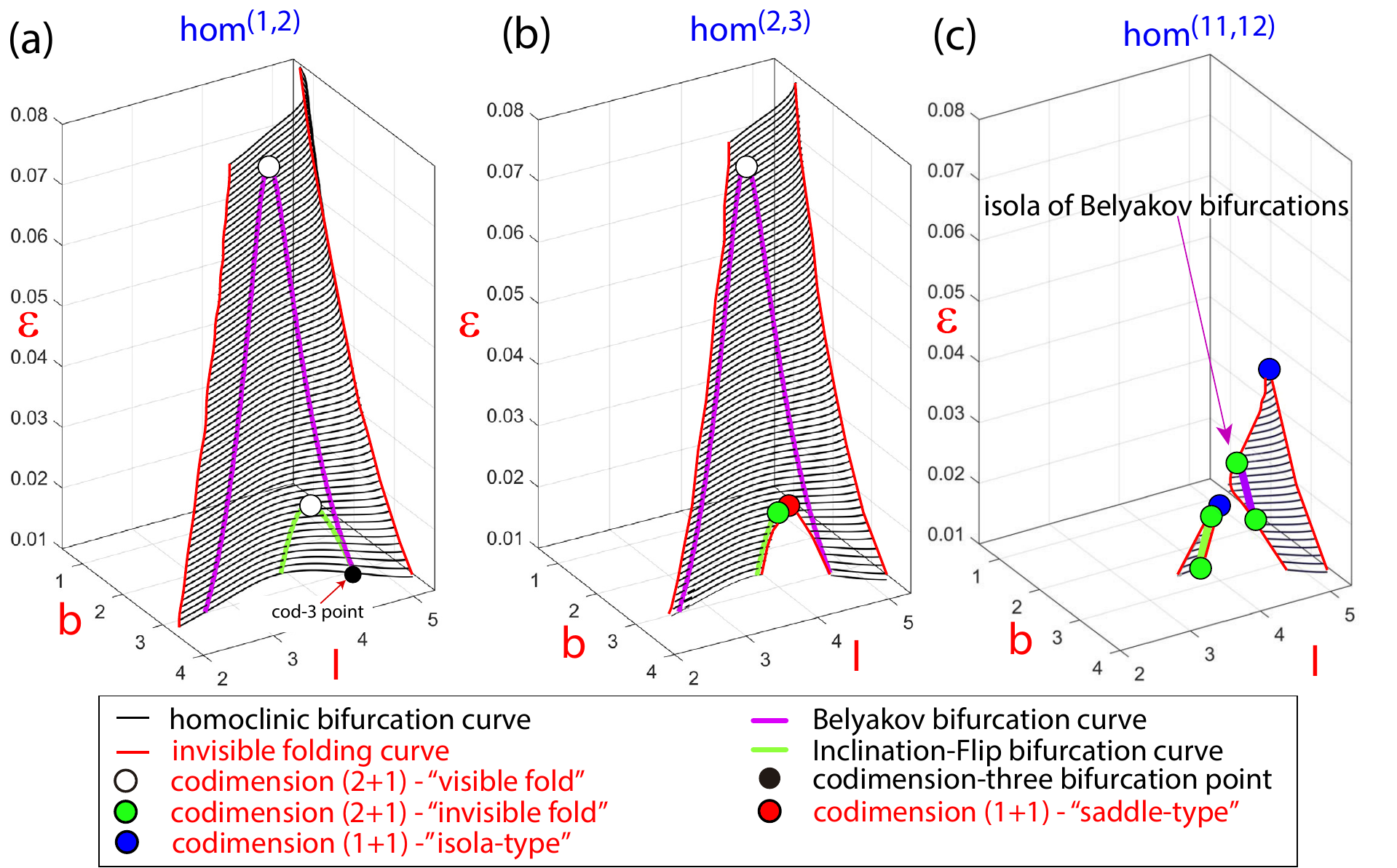}}
\caption{Three-parameter plot showing codimension-one homoclinic bifurcation surfaces. In plots (a), (b) and (c) it is presented one homoclinic bifurcation surface with different number of spikes, $hom^{(1,2)}$, $hom^{(2,3)}$ and $hom^{(11,12)}$, respectively. Some codimension-two homoclinic bifurcation curves are shown over the surfaces and some geometric bifurcations are pointed. }
\label{fig:1+1-extreme-homoclinic}
\end{figure*}

From the analysis in \cite{BIP2020}, it follows that there exist many codimension-one homoclinic bifurcation surfaces which are exponentially close each other and its number grows to infinity when the small parameter tends to zero. Moreover, since these homoclinic surfaces are tubular, the intersection of each surface with horizontal planes produce isolas (closed curves). These isolas exhibit a pair of extremely sharp folds and their width is also exponentially small (compare with the scheme given in Figure~\ref{topologia-def} with $d \ll 1$). Folding points determine two different sides in the isola and also in the bifurcation surface. Typically, for parameter values on one of the ``sides'' of the homoclinic bifurcation isola, the homoclinic orbit exhibits $n$ spikes and, for parameter values on the another face, $n+1$. This explains the notation $hom^{(n,n+1)}$ to refer the different isolas and surfaces. In Figure~\ref{fig:1+1-extreme-homoclinic} we show three-parameter plots showing some codimension-one homoclinic bifurcation surfaces and codimension-two homoclinic bifurcation curves computed using the continuation AUTO software (see \cite{AUTO2,AUTO}). To construct the surfaces, a collection of codimension-one homoclinic bifurcation curves are computed for different values of $\varepsilon$.  In  plots (a), (b) and (c) the $hom^{(1,2)}$, $hom^{(2,3)}$ and $hom^{(11,12)}$ surfaces are given, respectively, which include Belyakov and Inclination-Flip codimension-two bifurcation curves. This partial bifurcation diagrams illustrate, respectively, Cases I, II and III shown in the theoretical scheme provided in Figure \ref{topologia-def}.

Looking for codimension--one-plus-one geometric bifurcations, we pay attention to the bifurcation surfaces themselves. In Figure~\ref{fig:1+1-extreme-homoclinic}(c) we see how $hom^{(11,12)}$ splits into two disconnected components and each of them has a maximum (with respect to $\varepsilon$). These two points are isola-type codimension--one-plus-one geometric bifurcations with respect to $\varepsilon$. Also, as shown in plot (b) for the surface $hom^{(2,3)}$, a saddle-type codimension--one-plus-one geometric bifurcation is detected on the bifurcation surface. Most likely, isola-type codimension--one-plus-one geometric bifurcations are present in all homoclinic surfaces if the small parameter $\varepsilon$ grows enough and it is no longer ``small''. This fact explains why, as the small parameter grows, fewer bands of color appear in the 2D plots of Figure~\ref{fig:general} for the largest value of $\varepsilon$, indicating that the bursting orbits have fewer spikes.

On the other hand, codimension--two-plus-one geometric bifurcations correspond to folds, with respect to $\varepsilon$, that appear along curves of codimension-two homoclinic bifurcations and, as explained in the previous section, we distinguish between ``visible'' and ``invisible'' folds. Recall that in all two-dimensional manifolds of codimension-one homoclinic bifurcations we distinguish two leaves which are exponentially close (in fact the two leaves that form the isolas) and therefore they are indistinguishable in our visualization of the numerical results. They glue together along curves of sharp folding marked with a red line in Figure~\ref{fig:1+1-extreme-homoclinic}. When a curve of codimension-two homoclinic bifurcation folds inside one of the leaves, we get a \emph{codimension--two-plus-one visible fold}, but, if the fold is along one of the curves of folding of the whole surface (that is, going from one leaf to the another), and hence it is hidden for visualization, we have a \emph{codimension--two-plus-one invisible fold}. Note that visible folds appear in pairs, one on each leaf, but invisible folds are unique points. On the surface $hom^{(1,2)}$, both the Belyakov and the Inclination-Flip bifurcation curves show ``visible'' folds. On the surface $hom^{(2,3)}$,  the Belyakov bifurcation curve also undergoes a ``visible'' fold, but the Inclination-Flip curve presents now an ``invisible'' fold. And on the surface $hom^{(11,12)}$, both curves exhibit ``invisible'' folds.

The phenomenon of the ``invisible'' folds is a direct consequence of the existence of very thin tubular structures (and therefore isolas when a two-parameter section is considered) where the two leaves are infinitesimally close each other. Thus, if a codimension-two homoclinic bifurcation curve reaches the homoclinic surface folding curve, then it continues to the other side because the conditions in the phase space that are required to have the codimension-two bifurcation are still satisfied on the other side.

\begin{figure*}
\centerline{\includegraphics[width=1.\textwidth]{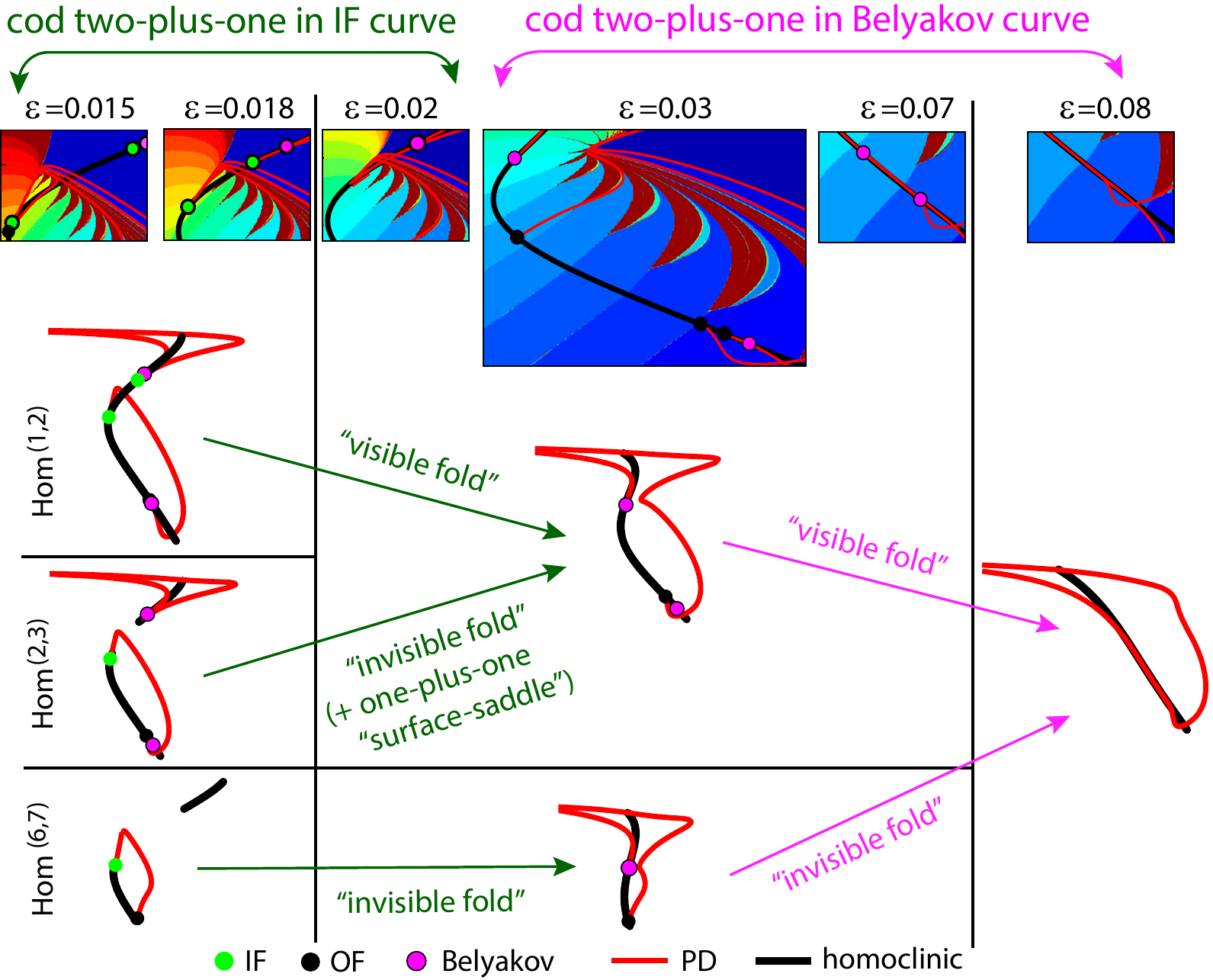}}
\caption{Schematic processes of some codimension--two-plus-one bifurcations in the HR model increasing the parameter  $\varepsilon$. On the top, SC pictures with some bifurcations illustrate the global panorama at several values of $\varepsilon$.
On the bottom pictures some numerical bifurcation curves provide an outline of different ways, through ``visible'' vs. ``invisible'' folds, to arrive to similar global situations.}
\label{fig:tabla}
\end{figure*}

In the Introduction we show in Figure~\ref{fig:general} how some codimension-two points disappear and in Section \ref{sec:01} we also explain that a subsidiary effect was the fact that different curves of periodic bifurcations may connect together, as illustrated in Figure~\ref{fig:general}(c).
Now we can connect all these phenomena with some geometric bifurcations. In Figure~\ref{fig:tabla} we show some codimension--two-plus-one bifurcations in the HR model increasing the small parameter  $\varepsilon$. We illustrate schematically the global structure before and after the codimension--two-plus-one bifurcations
and how the size of the homoclinic isolas determines the type: ``visible'' or ``invisible'' fold. As $\varepsilon$ increases, the first geometric bifurcation
in Figure~\ref{fig:general} is a
codimension--two-plus-one bifurcation that occurs on the Inclination-Flip (IF) homoclinic bifurcation curve. For the $hom^{(1,2)}$ case, as the homoclinic surface is big enough, there is a maximum of each of the IF curves (one on each leaf of the tubular surface), and so we have a ``visible'' fold and we observe the geometric ``collision'' of two pairs of IF points, one pair on each of the leaves of the homoclinic surface (the white points on the Figures~\ref{topologia-def} and \ref{fig:1+1-extreme-homoclinic}). On the contrary, for the rest of homoclinic surfaces, the IF bifurcation curves have ``invisible'' folds, as they are smaller and the surfaces are composed of isolas disconnected for small values of the parameter. Moreover, when $n$ grows enough $hom^{(n,n+1)}$ has only one branch of the IF curve, that is, IF points are only present on one of the two disconnected components of the isolas for a small value of $\varepsilon$. In these cases, the IF point seems to disappear on the limit of a homoclinic curve, and what really happens is that geometrically ``collides'' with the corresponding point of the another leaf of the surface that we cannot see (the green points on the Figures~\ref{topologia-def} and \ref{fig:1+1-extreme-homoclinic}). The evolution of the Belyakov curve is similar but for higher values of the small parameter $\varepsilon$.

Figure \ref{fig:tabla} also illustrates how codimension-one bifurcation curves arising from codimension-two bifurcation points are affected by the presence of a geometric bifurcation. For instance, following the surface $hom^{(1,2)}$, first a collision of two IF points is observed (a codimension-two-plus-one geometric bifurcation) which is accompanied  by the reconnection of two branches of period doubling bifurcation curves.

\begin{figure*}
\centerline{\includegraphics[width=1.\textwidth]{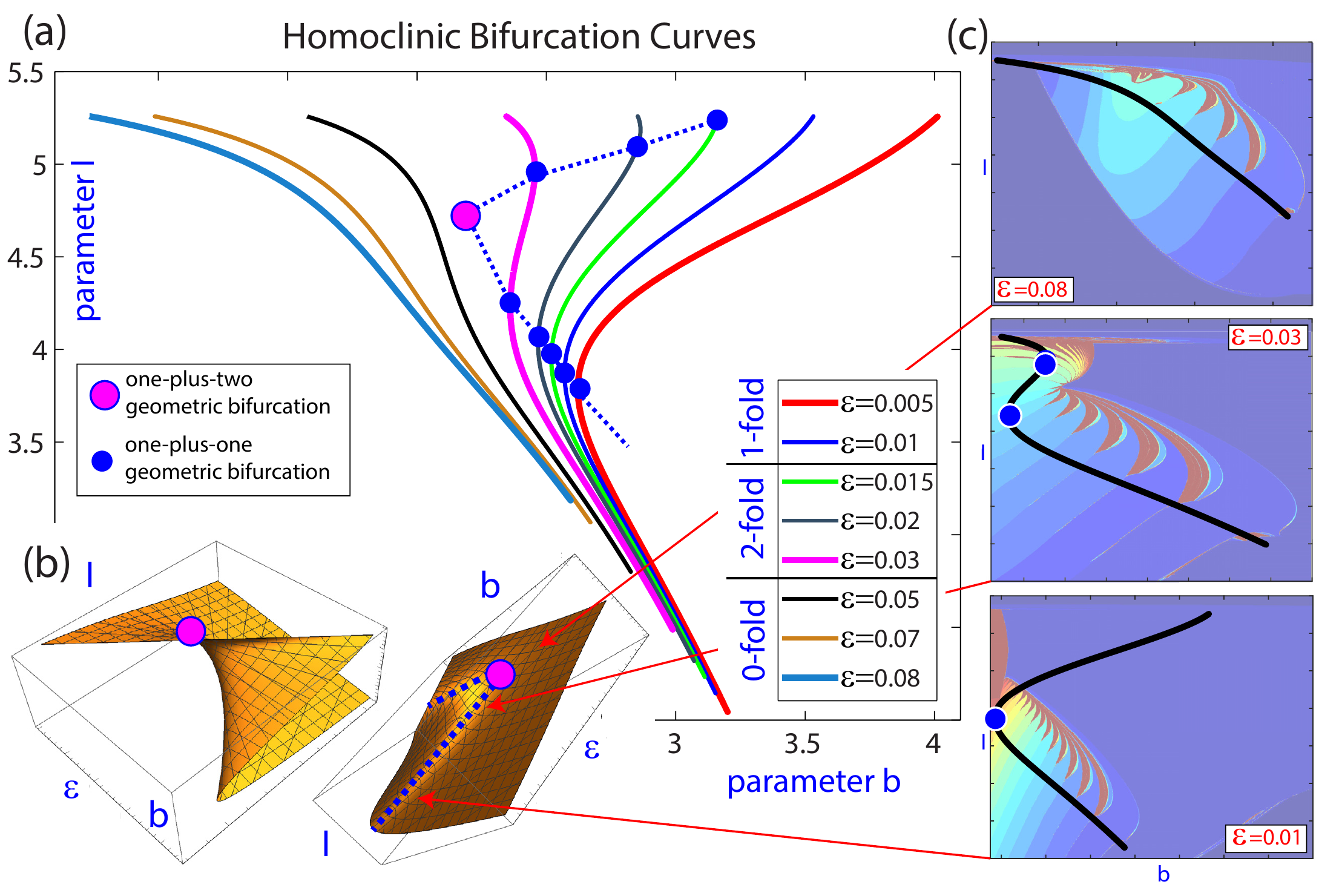}}
\caption{One-plus-one and one-plus-two homoclinic geometric bifurcations in the HR model. (a) $(b, I)$ projections of the $hom^{(1,2)}$ bifurcation curve for different values of $\varepsilon$. The location of the codimension--one-plus-one bifurcation points (``visible folds'') and the codimension--one-plus-two bifurcation point are given. (b) Schematic picture of the geometry of the $hom^{(1,2)}$ bifurcation surface presenting a \emph{cusp catastrophe} geometry giving rise to the  codimension--one-plus-two bifurcation point. (c) SC plots with the $hom^{(1,2)}$ bifurcation curves for three values of $\varepsilon$ showing the change of their geometry.}
\label{fig:oneplustwo}
\end{figure*}

Finally, to complete the panorama of geometric bifurcations given in Section~\ref{sec:1}, we show in Figure~\ref{fig:oneplustwo} how the HR model exhibits a codimension--one-plus-two geometric bifurcation point (compare with Figure~\ref{cod_1p2}, where the $\varepsilon$ parameter of HR is the parameter $\varepsilon_2$. Recall that, in spite of the notion of geometric codimension two has not been formally introduced, we have linked that concept to conceiving the parameter space as a two-parameter family of one-parameter bifurcation diagrams. So, when we describe this codimension--one-plus-two geometric bifurcation point, geometric bifurcations must be understood as tangencies or degenerated tangencies between straight lines in the parameter space (where the value of two parameters is fixed) and a bifurcation surface. On the plot (a) we show the projection of the $hom^{(1,2)}$ bifurcation curve for different values of $\varepsilon$ on the $(b, I)$ parameter plane. The location of the codimension--one-plus-one bifurcation points ("visible folds") on this projection plane is given and we see how they form a pair of geometric bifurcation curves. These curves have a cusp-type contact giving rise to a codimension--one-plus-two geometric bifurcation point. A scheme of the geometry of the $hom^{(1,2)}$ bifurcation surface is given in the two plots (b) presenting a \emph{cusp catastrophe} geometry giving rise to this codimension--one-plus-two bifurcation point (compare with the illustration provided in Figure~\ref{cod_1p2}). Note that as the homoclinic curves finish cutting this
codimension--one-plus-one curve, the cusp catastrophe surface is incomplete in some parametric regions giving just one fold point ("C" shape" vs. "Z" shape). Plots (c) show SC pictures with the $hom^{(1,2)}$ bifurcation curves for three values of $\varepsilon$ to see more clearly the process, and illustrating the changes in the geometry of the homoclinic curves.

\subsection{The FitzHugh-Nagumo system}
\label{ssB}
The FitzHugh-Nagumo system \cite{Fitz61,Nagumo} is a simplified version of the Hodgkin-Huxley model \cite{HH52} for the propagation of nerve impulses in axons, and the reduced ODE system can be written in the form~\cite{CKKOS2007}:
\begin{equation}
\label{FNmodel-ode}
\left\{\begin{array}{rcl}
U'&=&V,
\\
V'&=&\displaystyle{\frac{1}{\Delta}}\big(sV-U(U-1)(\alpha-U)+W-p\big),
\\[6.pt]
W'&=&\displaystyle{\frac{\varepsilon}{s}}(U-\gamma W).
\end{array}\right.
\end{equation}

System (\ref{FNmodel-ode}) has been extensively studied in the literature~\cite{CS2015,CS2018,CKKOS2007,guckue2009,guckue2010,KSS1997}. In \cite{CKKOS2007}, taking $s$ and $p$ as bifurcation parameters, it was shown the existence of $C$-shaped curves of homoclinic bifurcations (travelling waves in the original PDE FitzHugh-Nagumo system correspond to homoclinic orbits in system~(\ref{FNmodel-ode})). In fact it was observed that these $C$-shaped curves were isolas of an exponentially small wide, that authors called ``homoclinic bananas''. They also showed how ``bananas'' could split due to the existence of codimension-two bifurcation points on the homoclinic curve that played the role of terminal points for the branches following after the sharp turning points. That is, a ``banana split'' consists of two isolas joined by a curve. Geometrical explanations for the sharp turns were discussed in \cite{CS2015,guckue2009,guckue2010}. In \cite{CS2018}, the transition was analytically described using geometric singular perturbation theory and blow-up techniques.

We will show that the bifurcation diagram of system~(\ref{FNmodel-ode}) with respect to parameters $(\alpha,s,\varepsilon)$ exhibits geometric bifurcations which are related with the split of the homoclinic banana.

\begin{figure*}
\centerline{\includegraphics[width=1.\textwidth]{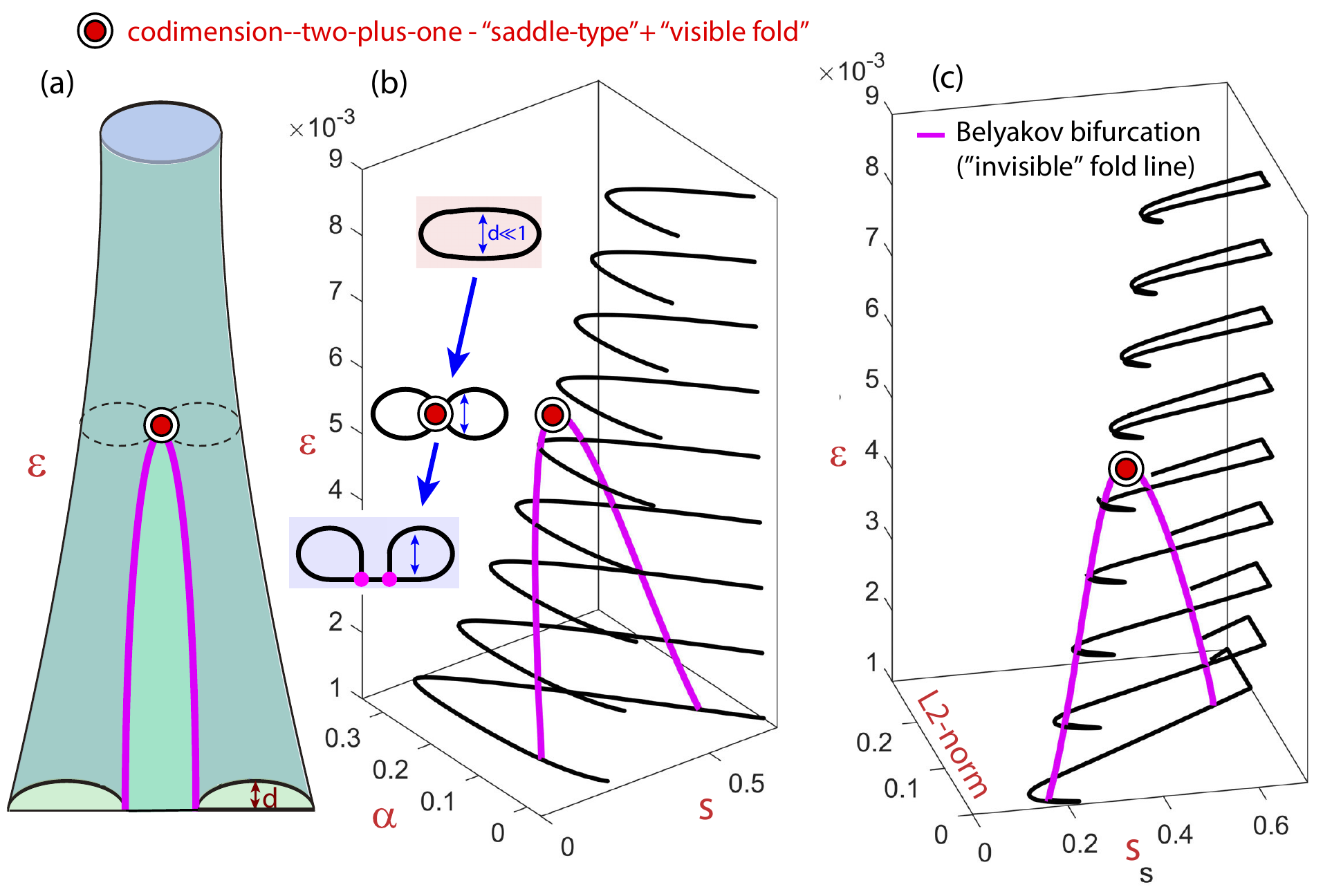}}
\caption{Three-parameter plot showing codimension-one homoclinic $hom^{(1,2)}$ bifurcation surfaces for the FitzHugh-Nagumo system. In plot (a) it is shown the theoretical structure of the homoclinic bifurcation surfaces, in this case it corresponds to the ``duck-foot shape'' case IIb of Figure~\ref{topologia-def}. Plots (b) and (c) show the homoclinic bifurcation curves for several values of the small parameter $\varepsilon$ and the codimension-two Belyakov bifurcation curve
obtained using the continuation software AUTO. Plot (b) is given in the three-parameter space $(\alpha, s, \varepsilon)$ while plot (c) is represented by using the parameters $(\alpha, \varepsilon)$ and the AUTO $L_2$-norm. }
\label{fig:FH}
\end{figure*}

In Figure~\ref{fig:FH} we present a three-parameter plot of the codimension-one homoclinic $hom^{(1,2)}$ bifurcation surfaces for the system~(\ref{FNmodel-ode}) with $p, \gamma=0$, $\Delta=1$, and considering $\alpha$, $s$ and $\varepsilon$ as bifurcation parameters. Plot (a) shows the theoretical homoclinic structure that corresponds to ---\emph{Case IIb}---, ``duck-foot surface'', which is obtained from the combination of a ``pair of pants'' and a ``disc'' (see Figure~\ref{topologia-def}). In this case a codimension-two Belyakov bifurcation curve acts as the limit bound of the tubular structures with a secondary homoclinic (after obtaining an extra spike) that returns to the Belyakov points and enables the continuation in a disc surface. That is, a transversal cut of the structure gives an isola for large values of $\varepsilon$;  for an intermediate value of $\varepsilon$, it gives a point of generation of two coupled isolas that coincides with the geometric maxima of the Belyakov curve (a codimension–two-plus-one bifurcation); and for small values of $\varepsilon$, it gives a set of two connected isolas, the ``homoclinic banana split'' that is described in the homoclinic literature. Plots (b) and (c) provide numerical
results obtained using the continuation software AUTO, and they show the homoclinic bifurcation curves for several values of the small parameter $\varepsilon$ and the codimension-two Belyakov bifurcation curve. Plot (b) is given in the three-parameter space $(\alpha, s, \varepsilon)$ and apparently we only observe just one curve for each $\varepsilon$, and no loop at all, but this is due to the very small distance between the sides of the isolas ($d \ll 1$). In order to observe the isolas and the connected isolas, we have to use the AUTO $L_2$-norm as it is shown in plot (c) using the parameters $(\alpha, \varepsilon)$ and the $L_2$-norm.

\section{Conclusions}
In this paper, partially bringing together previous studies, we introduce the concept of geometric bifurcation and illustrate how it appears in the prominent context of neural models. Of course, the Hindmarsh-Rose and FitzHugh-Nagumo systems are not unique examples. Geometric bifurcations are also expected to appear in other neural models and also in models coming from contexts other than Neuroscience (see examples in the already mentioned references \cite{algaba2003,algaba2011,algaba2015,Wieczorek2007,Wieczorek05}).

It should be noted that the study of geometric bifurcations may be of special interest in the framework of fast-slow systems. In this case, those parameters that modulate the slow dynamics are distinguished parameters. These parameters are the ones which one should fix to get 2D slices. Therefore, they specify a very concrete way of exploring the parameter space. In fact, it is possible to observe variations in the bifurcation diagrams when the distinguished parameters change, that is, signs of the presence of geometric bifurcations. In summary, if a multi-parameter space needs to be explored, there is no other option than to work with $k$-parameter slices with $k=1,2$ or, at most, $k=3$, and the different phenomena that we explain in this article could surely be present.

We formally introduce the notion of geometric bifurcation. In particular, we utilize a terminology (codimension $n$-plus-$m$) that emphasizes the distinction between bifurcations in families of dynamical systems and in ``families of families''.
Moreover, we studied the geometric bifurcations of codimension--one-plus-one and two-plus-one in three-parameter spaces. Additionally, we described a codimension--one-plus-two geometric bifurcation.
We were able to explain the appearance of different geometric bifurcations by the combined use of the Morse classification of 2D manifolds and bifurcation theory.

Finally, we conclude that this approach provides a nice tool to explain the different changes observed
in the phase plane when changing a parameter (like in real systems) as a mixture of dynamic and geometric bifurcations.

\section*{Acknowledgements}
RB has been supported by the Spanish Research projects PGC2018-096026-B-I00 and PID2021-122961NB-I00, and the European Regional Development Fund and Diputaci\'on General de Arag\'on (E24-17R and LMP124-18). SI and LP have been supported by the Spanish Research projects  MTM2017-87697-P and PID2020-113052GB-I00. LP has been partially supported by the Gobierno de Asturias project PA-18-PF-BP17-072.

\bibliographystyle{spmpsci}

\end{document}